\documentclass[11pt]{article}
\usepackage{amsfonts, epsfig, amsmath, amssymb, color,amsthm}
\usepackage{textcomp}
\usepackage{paralist}
\usepackage[normalem]{ulem}

\usepackage[english]{babel}
\newcommand{\RR}{\mathbb{R}}

\newcommand{\mM}{\mathcal{M}}

\newcommand{\vp}{\varphi}

\newcommand{\no}{\noindent}

\newcommand{\ra}{\rightarrow}

\theoremstyle{plain} 
\newtheorem{theorem}{Theorem}[section]
\newtheorem{corollary}{Corollary}[section] 
 
\newtheorem{proposition}{Proposition}[section] 

\theoremstyle{definition} 
\newtheorem{definition}{Definition}[section]

\newtheorem{condition}{Hypothesis}

\theoremstyle{remark} 
\newtheorem{remark}{Remark}[section]
\theoremstyle{definition}

\usepackage{lineno}

\DeclareMathSymbol{\ophi}{\mathalpha}{letters}{"1E}

\renewcommand{\phi}{\varphi}

\newcommand{\be}{\begin{equation}}
\newcommand{\ee}{\end{equation}}
\newcommand{\ben}{\begin{equation*}}
\newcommand{\een}{\end{equation*}}

\newcommand{\ba}{\begin{equation}\begin{aligned}}
\newcommand{\ea}{\end{aligned}\end{equation}}

\newfont{\cyrfnt}{wncyr10}
\def\J3{\cyrfnt{\rm \u{\cyrfnt I}}}
\def\j3{\cyrfnt{\rm \u{\cyrfnt i}}}

\usepackage[]{color}
\definecolor{DarkGreen}{rgb}{0.1,0.7,0.3}   

\definecolor{DarkGreen}{rgb}{0.1,0.7,0.3}   

\allowdisplaybreaks[4]

\begin{document}

\title{The Hele-Shaw free boundary limit 
\\ 
of Buckley-Leverett System }

\date{\null}

\author{Andr\'e 
Gomes$^1$, Wladimir Neves$^{2,\ast}$}

\date{}

\maketitle

\footnotetext[1]{ Instituto de Matem\'atica, Universidade Federal
do Rio de Janeiro, Brazil. E-mail: {\sl andredeoliveiragomes2@gmail.com.} }
\footnotetext[2]{ Instituto de Matem\'atica, Universidade Federal
do Rio de Janeiro, Brazil. E-mail: {\sl wladimir@im.ufrj.br.} \, $^\ast$ Corresponding author.}

\maketitle

 \begin{abstract} 
 This paper proposes a new approach to solving the Buckley-Leverett System,
 which is to consider a compressible approximation model characterized by a stiff 
 pressure law. Passing to the incompressible limit, 
the compressible model gives rise to a Hele-Shaw type free boundary limit 
of Buckley-Leverett System, and it is shown the existence of a weak solution of it.

 \end{abstract}

\noindent \textbf{Keywords:} 
\noindent \textbf{2010 Mathematical Subject Classification: } .
\tableofcontents
\newpage

\section{Introduction}

We are  interested
in this paper to study the Buckley-Leverett System
for multi-phase flows in porous media. 
This problem has been studied by many authors because, it plays an important role 
in a wide range of the planning and operation of oil wells, filtration problems, dam study, and 
also in a large class of interesting similar mathematical models 
(cf. \cite{BCBT, Chavent} and the references cited therein). 
The most prototypical example is the encroachment of water into an oil sand,
that is to say, the flow mixture of two phases (oil and water) in a porous medium, 
see Morris Muskat \cite{Muskat}.  

\medskip
The Buckley-Leverett System consists in
a conservation law that expresses the mass balance for the density $u(t,x)$
of the fluid-mixture evolving according to the nonlinear continuity equation
$$
 \partial_t u + \text{div}(g(u) \, \textbf{v})=0,
$$
where $g(u)$ represents an additional volume filling term in the density flux, and
$\textbf{v}(t,x)$ is the (seepage) velocity field given by 
$$
 h(u) \textbf{v} = - \nabla p,
$$
obtained from the empirical Darcy's law, (applied for each component of the mixture),
which describes the dynamics of the flow in relation 
with the scalar function $p(t,x)$, usually called pressure. Indeed, we observe that here 
$p(t,x)$ is not given by a law of state, but it is related to an internal constraint which is the incompressibility
condition for the velocity field of the flow, that is, 
$\text{div}(\textbf{v})= 0$.    
Then, we may write the Buckley-Leverett System as 
\begin{equation} 
\label{eq: the BL system}
\left \{
\begin{aligned}
&\partial_t u - \text{div}(\Theta(u) \nabla p)= 0,
\\[5pt]
& \text{div} \Big(\frac{\nabla p}{h(u)}\Big)= 0,
\end{aligned}
\right . 
\end{equation}
where $\Theta(u)= g(u) / h(u)$, with $h(u) \geq h_0$ for some $h_0> 0$. 
Here, $h^{-1}$ has the physical meaning of the total mobility, and 
we recall that \eqref{eq: the BL system} is also called the weak formulation of
multidimensional Muskat problem. 

\medskip
The system \eqref{eq: the BL system} is posed in $(0,T) \times \Omega$, 
for any real number $T> 0$ and $\Omega \subset \RR^d$ be an open domain, where
$d$ is the spatial dimension. Since we need some source 
either in the evolution equation $(\ref{eq: the BL system})_1$ or at the boundary of $\Omega$, that is, 
in orther to avoid non-trivial limit evolution, the domain 
$\Omega$ must be bounded or an exterior domain. We take the bounded one, which  
seems more physicaly correct for the applications. Thus, we prescribe some 
influx and outflux of the flow at the boundary of $\Omega$, denoted by $\Gamma$, that is to say, 
\begin{equation}
\label{BCINC}
\big(\Theta(u) \, \nabla p \big) \cdot \mathbf{n}= f \quad \text{on $\Gamma_T$},
\end{equation}
where $\mathbf{n}(r)$ denotes the unitary outer normal at $r \in \Gamma$, $\Gamma_T= (0,T) \times \Gamma$,
and $f$ is a given normal boundary flux, such that $|f(t,r)| \geq f_b$ for some $f_b> 0$. 
Also, the system \eqref{eq: the BL system} should be 
supplemented with an initial data for the density, that is, 
\begin{equation}
\label{ICBLS}
u(0)= u_0 \quad \text{in $\Omega$,} 
\end{equation}
where $u_0(x)$ is a nonnegative given function. 

\subsection{Difficulties and previous works.} 

One remarks that, the mathematical problem of proving global well-posedness of weak solutions 
for the Buckley-Leverett System is a mathematical challenge. In fact, the main question of 
showing existence of weak solutions remains open for many decades. 

\medskip
First, we observe that the velocity $\textbf{v}$ related to the Buckley-Leverett System
is expected to be a $L^2$ vector field.
Therefore, we have to treat with a scalar conservation law, that is the equation $(\ref{eq: the BL system})_1$, 
in the class of roughly coefficients in order to solve system \eqref{eq: the BL system}. 
Moreover, due to the presence of discontinities (shock waves) in the solution $u$,
even if $u_0$ is a smooth function, the entropy formulation is usually adopted. Recall that, the 
entropy formulation is an admissible criteria to establish uniqueness of weak solutions, not a necessary 
condition to show existence of them. 

\medskip
Albeit, if one tries to use the kinetic formulation of scalar conservation laws
to show existence of weak solutions to the Buckley-Leverett System, then
the entropy formulation is needed. But, it is well known that we do not have 
uniqueness of weak solutions 
for the associated kinetic-transport equation 
when the drift vector field is just in $L^2$, 
(see \cite{{Depauw}, OKWN}).
Hence the mapping to apply a fixed point argument to resolve the
system \eqref{eq: the BL system} is not well defined. 
Moreover, from the lack of uniqueness 
the DiPerna--Lions, Ambrosio renormalization
property does not work, which is an important issue to show 
compactness of approaching sequences of regular solutions 
using, for instance, the rigidity result of Perthame-Dalibard 
\cite{perth-dali-tams}.

\medskip
We also observe that, the flux function $g(u) \textbf{v}$ is degenerated, that is, 
there exists $\xi \in \mathbf{S}^{d-1}$, such that the mapping 
$\lambda \mapsto g(\lambda) \textbf{v}(t,x) \cdot \xi$ is constant on non-degenerated intervals,
for $(t,x)$ in an open set contained in $(0,T) \times \Omega$. Therefore, there is no hope to 
apply the H-measures theory associated with the measure value solutions for conservation laws, 
see Panov \cite{PaARMA}. 
Similarly, it is not possible to consider the kinetic theory introduced by 
Lions, Phertame, Tadmor \cite{LPT2}. Indeed, 
a slightly stronger notion of non-degeneracy hypothesis than the previous one does not work either, that is, 
$$
   \text{$\exists (\xi_0,\xi) \in \mathbf{S}^d$, such that $\mathcal{L}^1(\{\lambda \in \RR ; \xi_0 + \xi \cdot g(\lambda) \textbf{v}(t,x) = 0\})> 0$, for $(t,x)$ in an open set}. 
$$

\medskip
To overcome the difficulties mentioned above,
the Buckley-Leverett System
has undergone significantly modifications
in some directions.
A quite general case can be 
described by the following partial differential equations: 
given non-negative (viscosity) functions $\mu$ and $\nu$, we have 
\begin{equation}
\label{eq.general}
\left \{
\begin{aligned}
    &\partial_t u + \text{div} \big(g(u) \,  \mathbf{v} \big)=  \text{div}(\mu \, \nabla u) + F(\cdot,u),
\\[5pt]
    &  - \text{div}( \nu \, \mathbf{D} \mathbf{v}) + h(u) \, \mathbf{v}= -\nabla p + \mathbf{G}(\cdot, u), 
\\[5pt]    
    & \text{div}(\textbf{v})= 0,
\end{aligned}
\right.
\end{equation}
where $F$ and $\mathbf{G}$ are external forcing terms, $\mathbf{D}(\cdot)$ is the symmetric part of the gradient operator,  
and the above system must be supplemented 
with suitable initial and boundary conditions. 
There exists a list of 
important papers retated to \eqref{eq.general}, in particular following the two main 
directions, that is, $\mu>0$, $\nu= 0$, and $\mu= 0$, $\nu> 0$. 

\medskip
In \cite{KruskSukor} Kru\u{z}kov, Sukorjanski\u{i} showed that, the system \eqref{eq.general} 
is well-posed in the class of semiclassical solutions for $\mu= \mu(x,u)> 0$, and $\nu= 0$, 
see Definition 2 in that paper. Hence they considered a parabolic--elliptic structure. 
Also in this direction, we mention the papers: 
Eymard, Herbin, Michel  \cite{Eymard_Herbin_Michel}, 
Benedetto, Gianazza, Vespri \cite{Benedetto_Gianazza_Vespri}, 
Michel \cite{Michel}, Lenzinger, Schweizer \cite{len}. 

\medskip
Conversely, Chemetov, Neves showed in \cite{ChemetovNeves2} 
existence of weak solutions for the system \eqref{eq.general} 
considering $\mu= F= \mathbf{G}= 0$, and
$\nu= const> 0$, see Definition 6.1 in that paper.
Thus a hyperbolic--elliptic structure was considered. 
Along the same direction, we mention the papers:
Colcite, Karlsen, Mishra, Risebro \cite{Colcite1}, and 
Colcite, Mishra, Risebro, Weber \cite{Colcite2}.  

\medskip
One remarks that, the preserved hyperbolic structure in the system \eqref{eq.general}  
is a noteworthy aspect, especially when sharp interfaces between the components of the misture is a crucial 
issue for applications. Moreover, it is important to note that
the provided list of references mentioned above is not exhaustive, 
and we address the reader for many further references cited therein.

\subsection{Motivation and strategies.}

In the present paper, we consider a new aproach to solve 
the Buckley-Leverett System \eqref{eq: the BL system},
which is the stiff approximation (or stiff limit solutions), and 
in particular, we are not going to apply the entropy formulation.
Our main motivation 
to use this approach comes from the paper
Perthame, Quirós, Vázquez \cite{PerthameARMA}. In that paper, they considered 
both incompressible and compressible  
mechanical models of tumor growth, where 
these two classes of models are related by using a stiff pressure law, 
see equation \eqref{stiffpressure} below. 
Under some hypotheses, they succeed to pass to the limit as $\gamma \to \infty$ (called incompressibility limit, or stiff pressure limit) 
in the compressible model, and obtain a weak solution to the incompressible model, which is
a free boundary problem of Hely-Shaw type. 

\medskip
Therefore, let us present our strategy here. 
We consider a (companion) compressible 
model to \eqref{eq: the BL system}, due to a convenient pressure equation of state.  
More precisely, 
given $\gamma>1$  let us consider the constitutive pressure having a power law in relation 
with the density of the fluid, that is to say, 
\begin{equation}
\label{stiffpressure}
p_\gamma = p_\gamma(u):= u^\gamma.
\end{equation}
Then, for $\alpha> 0$ (to be chosen a posteriori) we consider the following partial differential equation
\begin{equation} 
\label{eq: the compressible BL system}
\partial_t u_\gamma + \text{div}(\Theta(u_\gamma) \nabla p_\gamma)= \frac{1}{\gamma^\alpha} \, u_\gamma \, \Phi(p_\gamma)
\quad \text{in $(0,T) \times \Omega$}, 
\end{equation}
 supplemented with the boundary and initial conditions, given respectively by  
 \begin{equation}
 \label{BCIC}
\big(\Theta(u_\gamma) \, \nabla p_\gamma \big) \cdot \mathbf{n}= f_\gamma \quad \text{on $\Gamma_T$},
\quad \text{and} \quad u_\gamma(0)= u_{0\gamma} \quad \text{in $\Omega$},  
\end{equation}
where $f_\gamma$ and $u_{0\gamma}$ are precisely described in Hypothesis B of Section \ref{NotHyp}. 

\medskip
One observes that, 
the equation \eqref{eq: the compressible BL system} differs from the first equation 
of \eqref{eq: the BL system} by the introduction of an artificial forcing term rescaled by 
$\frac{1}{\gamma^\alpha}$, given by a smooth function $\Phi$, such that $\Phi' < 0$ and $\Phi(p_M)= 0$ 
for some threshold value $p_M> 0$ known in the literature as homeostatic pressure. The introduction 
of this term is inspired on the work \cite{PerthameARMA}, and despite we have introduced this forced term,
we are going to call \eqref{stiffpressure}, \eqref{eq: the compressible BL system}, and \eqref{BCIC}: compressible Buckley-Leverett System. 

 \medskip
 Under good assumptions on the initial data $u_\gamma(0)$ and also $f_\gamma$ in \eqref{BCIC}, 
 we expect that the family of solutions $\{(u_\gamma, p_\gamma)\}$ of the system \eqref{stiffpressure}, \eqref{eq: the compressible BL system} and \eqref{BCIC} will converge 
 when $\gamma \to \infty$, (in some suitable sense to be defined later), to an element $(u_\infty, p_\infty)$ solving \eqref{eq: the BL system}, \eqref{BCINC} in the weak sense, which is
 similarly to \cite{PerthameARMA}, a free boundary problem of Hele–Shaw type. 
 Moreover, the initial condition $u(0)= u_0$ in $\Omega$
is obtained in the $L^1-$strong sense. 

\medskip
Along the same strategy we address the reader to the following papers.
In \cite{Gil_Quiros} Gil, Quirós applyied the incompressibility limit to study the convergence (as $m \to \infty$) of solutions to the initial-boundary value problems 
for the pourous media equation, that is to say
$$
\partial_t u - \Delta u^m= 0 \quad \text{in $(0,\infty) \times \Omega$}
$$
to Hele-Shaw flow, 
where $\Omega$ is an open subset of $\RR^d$, bounded or an exterior domain. 
Similarly, the stiff pressure limit was also considered by Kim, Po\u{z}\'ar 
\cite{Kim}. In that paper, they study the assymptotic limit (as $m \to \infty$) 
of solutions to the Cauchy problem for the pourous media equation, 
that is
$$
\partial_t u - {\rm div}(u \nabla p)= u \, G(p) \quad \text{in $(0,\infty) \times \RR^d$}
$$
to Hele-Shaw flow, where $p= \frac{m}{m-1} u^{m-1}$ and $G$ has the same behavior as $\Phi$ described above.  
Recently, following the same strategies as cited above, 
the authors Guillen, Kim, Mellet in \cite{guillen2022hele} considered the porous media equation with a (nonmonotone) source term, that is
$$
\partial_t u - {\rm div}(u \nabla p)= \lambda \, u  \quad \text{in $(0,\infty) \times \Omega$},
$$
where $\Omega \subset \RR^d$ is an exterior domain and $\lambda(t,x)$ is a bounded function taking both positive or negative values.

\medskip
Moreover, the idea of stiff limit solutions (the passage from the compressible model to the solution of incompressible one)
was used in some other important papers. In \cite{Noemi_Benoit}  David, Perthame extended the so-called 
complementarity condition established in \cite{PerthameARMA} for purely mechanical model, taking into 
account the consumption of nutrients. Merllet, Perthame, Quirós \cite{Merllet_Perthame_Quiros} consider
weak solutions to a  problem modeling tumor growth, in particular, they show 
the regularity properties both for the solution and its free boundary. Finally, the Hele-Shaw type free boundary problem 
for a tumor growth was also studied by Perthame, Quirós, Tang, Vauchlet \cite{Vauchlet}. The novelty in that paper was the active 
of cells due to a difusion term, which gives more regularity to the density of cells, and also the free boundary velocity is no more
given by the gradient of the pressure.

\subsection{Notation and main hypothesis.}
\label{NotHyp}

Let $U$ be an open set in $\RR^d$.
We denote by $\mathcal{L}^d$, (or $dx$, $d\xi$, etc.),
the $d-$dimensional Lebesgue measure, and by $\mathcal{H}^\theta$, (or $dr$), the $\theta-$dimensional Hausdorff
measure. By $L^p(U)$ we denote the set of real $p-$summable functions 
with respect to the Lebesgue measure,
and the vector counterparts of these spaces are denoted by
$\mathbf{L}^p(U)= \big(L^p(U)\big)^d$. Similarly we have $L^p(\partial U; \mathcal{H}^{d-1})$, 
where $\partial U$ is the boundary of the set $U$. 

\medskip
The fractional Sobolev space is denoted by $W^{s,p}(U)$, where a 
real $s\geqslant 0$ is the
smoothness index, and a real $p\geqslant 1$
is the integrability index.
The notation $BV(U)$, $\mM^1(U)$ stands, respectively, for the space of bounded variation functions 
and the space of all finite Radon measures, defined on $U$.

\medskip
Definitely, we denote by $\Omega$ a bounded open set in $\RR^d$, with $C^2-$boundary denoted by $\Gamma$. 
For $T> 0$, we use the notations $Q_T= (0,T) \times \Omega$ for the time-space cylinder, 
and also $\Gamma_T= (0,T) \times \Gamma$ is the lateral boundary. 
We frequently use the abbreviated form $u:= u(t,x)$, $p:=p(t,x)$ for $(t,x) \in Q_T$. 
We also define the positive and negative part of a function $w$ as follows, $|w|^{\pm}:= \max \{ \pm w,0 \}$.
Moreover, we have 
$$
\operatorname{sgn}^{+}(w):=\left\{\begin{array}{ll}
1, & \text { if } w>0 \\
0, & \text { if } w \leqslant 0,
\end{array} \quad \operatorname{sgn}^{-}(w):=\left\{\begin{array}{ll}
\!\! -1, & \text { if } w<0 \\
\; 0, & \text { if } w \geqslant 0.
\end{array}\right.\right.
$$

\medskip
Now, we consider the main hypothesis of this paper. 
\begin{condition} \label{condition: coefficients of the equations} 
\no We assume the following hypothesis on the nonlinearity coefficients:
$$
\text{$g, h \in W^{1,\infty}_{\rm{loc}}(\RR)$, $h \geq h_0 >0$, $g(z)> 0$ for each $z> 0$, $g(0)= 0$, and}  
$$
\begin{equation}
\label{CondGPrime}
\text{for $\delta> 0$ sufficiently small, $|z-1| < \delta  \Rightarrow g^\prime(z) \geq 0$.}  
\end{equation}

Moreover, we assume that the following function 
\begin{equation}
\label{ThetaLambda}
\text{$\Theta(z):= \frac{g(z)}{h(z)}$ is non-decreasing,}
\end{equation}
which is bounded and Lipschitz on compact sets.  
\end{condition}

\begin{condition} \label{condition: initial condition for the compressible BL}
We assume the following hypothesis on the data: 

\smallskip
1. For some non-negative intial condition $u_0 \in L^1(\Omega)$, the 
initial data $u_{0 \gamma}$ satisfies
\begin{equation}
\label{DATA}
\left \{
\begin{aligned}
&u_{0\gamma} \geq 0,  \quad p_\gamma(0) \equiv p_\gamma(u_{0 \gamma}) \leq p_M, 
\\[5pt]
&\|u_{0 \gamma} - u_0 \|_{L^1(\Omega)} \ra 0 \quad \text{as $\gamma \ra \infty$}, 
\end{aligned}
\right . 
\end{equation} 
and further 
\begin{equation}
\label{BONDGRADU0}
\sup_\gamma ||\nabla u_{0 \gamma}||_{L^1(\Omega)} < \infty.
\end{equation}

2. For some $f \in L^\infty(\Gamma_T; \mathcal{H}^d)$, there exists $f_\gamma$ defined in $\RR^{d+1}_+=(0,\infty) \times \RR^d$ such that 
\begin{equation}
\label{DATABOUNDARY}
\|f_{\gamma} - f \|_{L^1(\Gamma_T; \mathcal{H}^d)} \ra 0 \quad \text{as $\gamma \ra \infty$},
\end{equation}
and also for any compact set $K \subset \RR^{d+1}_+$, 
\begin{equation}
\label{BONDFLUX}
\sup_\gamma ||\nabla f_{\gamma}||_{L^1(K)} < \infty.
\end{equation}
Moreover, we assume without loss of generality that, $u_{0\gamma}, f_\gamma$ are smooth functions. 
 \end{condition}

 \begin{remark}
 One observes that {\bf Hypothesis A} is quite general. In particular, the function $\Theta(\cdot)$ be non-decreasing 
 is a usual assumption, see for instance  \cite{Eymard_Herbin_Michel}. 

Moreover, the  {\bf Hypothesis B} is satisfied by a large set of functions. Indeed, $BV \cap L^\infty$ is contained in it. 
\end{remark} 

\section{The Statement of Main Result}

In what follows we define what is meant to be a weak solution of the Buckley-Leverett System. 
\begin{definition}
\label{DEFSOL} 
A pair of functions $(u,p)$, with $u \in L^\infty(Q_T)$ and  $\nabla p \in L^{2}(Q_T)$, 
is called a weak solution to the the Buckley-Leverett System \eqref{eq: the BL system}, 
if it satisfies
\begin{align}
\iint_{Q_T}& \big(u \;\phi _{t} + \,\Theta(u) \; \nabla p \cdot \nabla \phi \big)\ d{x} \,dt = \int_{\Gamma_T} f \, \phi \, d\mathcal{H}^d,  \label{DGV21}
\end{align}
for any test function $\phi \in C_{c}^{\infty}((0,T)\times \mathbb{R}^{d})$, and also for a.
a. $t \in (0,T)$ 
\begin{equation}
\int_{\Omega} \frac{\nabla p(t)}{h(u(t))} \cdot  \nabla \xi \, dx= 0, \label{V21}
\end{equation}%
holds for any $\xi \in C^\infty_c(\Omega)$. 
Moreover, the initial data for the density is attained in the $L^1-$strong sense, that is 
$$
{\rm ess}\!\! \lim_{t \to 0^+} \int_{\Omega} |u(t) - u_0| \, dx= 0. 
$$
\end{definition}

\medskip
The main result of the paper is the following
\begin{theorem} [Main Theorem] 
\label{thm: main theorem}
Fix $\alpha >1$.  
Under Hypothesis \ref{condition: coefficients of the equations}--\ref{condition: initial condition for the compressible BL}, it follows that:
\begin{enumerate}
\item The unique strong solution $(u_\gamma, p_\gamma)$ of the system \eqref{stiffpressure}, \eqref{eq: the compressible BL system} and \eqref{BCIC}
converges strongly in $L^1(Q_T)$ 
as $\gamma \ra \infty$ to an element $(u_\infty, p_\infty)$, with $u_\infty \in C([0,T]; L^1(\Omega)) \cap BV(Q_T)$ and $p_\infty \in BV(Q_T)$.
\item The functions $u_\infty$ and $p_\infty$ are such that $0 \leq u_\infty \leq 1$ and $0 \leq p_\infty \leq p_M$ pointwise for almost all 
$(t,x)$ in $Q_T$.  
Also $u_\infty$ and $p_\infty$ satisfy the graph relation $p_\infty (u_\infty -1) =0$.

\item The following weak convergence in $L^2(Q_T)$ holds:
\begin{align*}
\nabla p_\gamma \rightharpoonup \nabla p_\infty \quad \text{ as } \gamma \ra \infty.
\end{align*} 
\item The pair $(u_\infty, p_\infty)$ solves \eqref{eq: the BL system}, \eqref{BCINC} in the sense of Definition \ref{DEFSOL} with initial data $u_0$.
\item The further important qualitative properties for the limit holds
$$
   \partial_t u_\infty, \partial_t p_\infty \geq 0 \quad \text{in $\mathcal{D}^\prime(Q_T)$}, 
$$ 
that is, the density is on the rise, accompanied by a corresponding increase in pressure.
\end{enumerate}
\end{theorem}

\begin{remark}
We do not discuss in this paper the important issue of regularity of the weak solutions, and also the free boundary.
Neither, it is mentioned something concerning the uniqueness. These two subjects are also of our main interest. 
\end{remark} 

\subsection{Strategy of the proof. }

As previously communicated, we approximate the Buckley-Leverett System \eqref{eq: the BL system}, \eqref{BCINC}, and \eqref{ICBLS}
by a companion system given by \eqref{stiffpressure}, \eqref{eq: the compressible BL system}, and \eqref{BCIC}. The introduction 
of the source term $u_\gamma \, \Phi(p_\gamma)$ at the scale $1/\gamma^\alpha$ (with $\alpha>1$) facilitated us to use techniques of passing to the stiff limit as 
$\gamma \ra \infty$ inspired by the work \cite{PerthameARMA}, (see also \cite{Noemi_Benoit}),  and methods such as Hele-Shaw limits for porous medium equations, e.g. 
\cite{Gil_Quiros, guillen2022hele, Kim}.
 
 \medskip
Now, the literature in generalized porous media equations, that is, the filtration equation with a non-linear flux function, (see \eqref{eq: the system for comparison} below), 
assures the existence of the strong solution for \eqref{stiffpressure}, \eqref{eq: the compressible BL system} and \eqref{BCIC}.  
This family of strong solutions will converge in the stiff limit for a weak solution of \eqref{eq: the BL system} in the distributional sense. 
Indeed, to pass to the limit in \eqref{eq: the compressible BL system} as $\gamma \ra \infty$ and recover the weak solution of 
$$
    \partial_t u - \text{div} (\Theta(u) \nabla p)=0,
$$ 
one may derive uniform bounds in $W^{1,1}(Q_T)$ for $\{u_\gamma\}_{\gamma >1}$ and $\{p_\gamma\}_{\gamma >1}$.  
The fact that we are working on a bounded domain with $C^2-$regular boundary 
and due to the compact imbedding of $W^{1,1}$ in $L^1$, it follows that 
$u_\gamma \ra u_{\infty}$ and $p_\gamma \ra p_{\infty}$ strongly in $L^1(Q_T)$ as $\gamma \ra \infty$. Moreover, 
we have $u_{\infty}, p_\infty \in BV(Q_T)$. 

\medskip
To conclude such uniform bounds with respect to $\gamma >1$ we derive a crucial technical estimate of Benilan-Aronson type,
see Section \ref{Benilan-Aronson type estimates}. It is interesting to notice that the nonlinear perturbation in the divergence of the 
gradient of the pressure as oposed to the linear perturbation considered in \cite{PerthameARMA} makes the derivation 
of our Benilan-Aronson inequality much more intrincate.  This estimate plays a crucial role in the derivation of the 
uniform bounds in $\gamma>1$ for the spatial gradients of the density $u_\gamma$ and the pressure $p_\gamma$.  
This estimate reveals again its power when allowing to prove the weak compactness of the gradient of the pressure 
$\nabla p_\gamma$ in $L^2(Q_T)$.  Hence the uniqueness of the limit in $BV$ allow us to recover that, the limiting 
$\nabla p_\infty$ is in $L^2(Q_T)$.  

\medskip
Finally we would like to remark that the incompressibility condition 
$$
   \text{div} \Big ( \frac{\nabla p}{h(u)}\Big )=0
$$ 
is stated almost surely in $[0,T]$ and in the weak sense in $\Omega$,  i.e. multiplying by test functions in $\Omega$.  
It is in its proof where the Hele-Shaw limit,  that is obtained for the limiting density $u_\infty$ and limiting pressure 
$p_\infty$, plays a pivotal role in streamlining the proof  
and enabling the use of energy estimates for the equation of 
the perturbed pressure $p_\gamma$.  We would like to emphasize that the incompressibility condition of \eqref{eq: the BL system}
should not be confused with the complementarity condition obtained in the works \cite{Noemi_Benoit, PerthameARMA} and other references.  
The mathematical formulations are different and evidently the physical phenomena that they model are distinct.

\section{Elements for the Proof of Main Theorem} 

In this section, we derive the important results to establish the proof of Theorem \ref{thm: main theorem}. 
To this end, we mainly consider the compressible Buckley-Leverett System. 

\subsection{Existence and uniqueness of solutions for compressible model.}

This section aims to provide the existence, uniqueness and regularity needed for the function $u_\gamma$ that solves 
the system \eqref{stiffpressure}, \eqref{eq: the compressible BL system} and \eqref{BCIC}, and enable us to obtain all the uniform estimates needed in the sections that follow.

\begin{theorem}
\label{thm: existence uniq strong solution of the compressible BL}
Let $\alpha>1$ be fixed arbitrarily, and let $\gamma>1$. Under 
Hypothesis \eqref{condition: coefficients of the equations}--\eqref{condition: initial condition for the compressible BL},
there exists a unique pair $(u_\gamma, p_\gamma)$, 
solving the system \eqref{stiffpressure}, \eqref{eq: the compressible BL system} and \eqref{BCIC}
with initial-condition $u_{0\gamma}$ and normal boundary flux $f_\gamma$. 
Moreover, one has for every $\gamma>1$ that
\begin{align} \label{eq: comparison}
0 \leq u_\gamma \leq \sqrt[\gamma]{p_M} \quad \text{ and } \quad 0 \leq p_\gamma \leq p_M \quad \text{a.e.}-Q_T.
\end{align}
\end{theorem}

\begin{proof}
First, let $(u_{0\gamma}, f_\gamma)$ be satisfying Hypothesis \ref{condition: initial condition for the compressible BL}. 
Then, we consider the system  \eqref{stiffpressure}, \eqref{eq: the compressible BL system} and \eqref{BCIC} written in the following form
\begin{equation}
\label{eq: the system for comparison}
\left \{
\begin{aligned}
\partial_t u_\gamma= \Delta \Psi(u_\gamma) &+ \frac{1}{\gamma^\alpha} u_\gamma \Phi(p_\gamma), \quad (p_\gamma= u_\gamma^\gamma), 
\\
\nabla \Psi(u_\gamma) \cdot \mathbf{n}&= f_\gamma, 
\\
u_\gamma(0)&= u_{0\gamma},
\end{aligned}
\right .
\end{equation}
where the function $\Psi(\cdot)$ is defined by 
\begin{equation}
\label{PSIGH}
\Psi(z):=  \gamma  \int_0^z  \Theta(s) \, s^{\gamma-1} ds.
\end{equation}

\medskip
One observes that, the equation $(\ref{eq: the system for comparison})_1$ is a generalized porous media equation,  
and the function $\Psi$ is a $C^{1,1}-$continuous increasing function and superlinear at infinity. 
Indeed, the possibly degenerate case follows from replacing $\Psi(z)$ with 
$\Psi_\gamma(z)= \Psi(z) + \frac{1}{\gamma^\alpha} z$. 
Morevoer, we recall that $\Psi^\prime$ is differentiable a.e. by applying Rademacher's Theorem. 
Also we have $\nabla \Psi(u_\gamma) \cdot \mathbf{n} \equiv  \Theta(u_\gamma) \nabla p_\gamma \cdot \mathbf{n}$. 
Then, from the theory developed by Ph. Benilan in \cite{ph1972equations} 
and lately by Ph. Benilan, M. G. Crandall, P. Sacks \cite{benilan1988some} (see also 
Vázquez \cite{Vazquezbook},  Chapter 11),    
we infer that there exists a unique pair $$(u_\gamma, p_\gamma) \in  C([0,T]; L^1(\Omega)) \times C([0,T]; L^1(\Omega))$$
solving weakly the initial boundary value problem \eqref{eq: the system for comparison}. 

\medskip
Finally, under such assumptions and since 
the equation $(\ref{eq: the system for comparison})_1$ is parabolic non-degenerate 
we can apply the standard (quasilinear) parabolic theory and obtain existence and uniqueness of strong
solutions for the problem \eqref{eq: the system for comparison}. Moreover, we obtain \eqref{eq: comparison}
from comparison principle, that is, the maximum principle in the nonlinear context. 
\end{proof}

\subsection{$L^1$ bounds for $u_\gamma$ and $p_\gamma$.}

\begin{proposition} 
\label{prop_ L1 bounds pgamma, ugamma} Let us assume that the set of Hypotheses \ref{condition: coefficients of the equations}-\ref{condition: 
initial condition for the compressible BL} hold.  Fix $\alpha>1$ and for every $\gamma>1$ let $u_\gamma$ 
be the unique solution of the system \eqref{stiffpressure}, \eqref{eq: the compressible BL system} and \eqref{BCIC} in the sense of Theorem \ref{thm: existence uniq strong solution of the compressible BL}. 
Then, for any $t \in [0,T]$, 
\begin{align} \label{eq: L1 estimate ugamma}
||u_\gamma(t)||_{L^1(\Omega)} \leq e^{\frac{\Phi(0)T}{\gamma^\alpha}} \big( ||u_{0\gamma}||_{L^1(\Omega)}+ 2  \int_{\Gamma_T} |f_\gamma| \ drdt \big), 
\end{align}
and 
\begin{align} \label{eq: L1 estimate pgamma}
||p_\gamma(t)||_{L^1(\Omega)} \leq (p_M)^{(\gamma-1)/\gamma} e^{\frac{\Phi(0)T}{\gamma^\alpha}} \big( ||u_{0\gamma}||_{L^1(\Omega)} 
+ 2  \int_{\Gamma_T} |f_\gamma| \ drdt \big).
\end{align}
\end{proposition}

\begin{proof}
1. First, let $u_\gamma$ be the unique solution of \eqref{eq: the compressible BL system}, that is to say, 
\begin{equation}
\label{EQUGAMMA}
\begin{aligned}
\partial_t u_\gamma &= \text{div} \Big ( \Theta(u_\gamma)  \nabla p_\gamma \Big ) + \frac{1}{\gamma^\alpha} u_\gamma \Phi(p_\gamma)
\\[5pt]
&= \gamma \, \text{div} (\Theta(u_\gamma) u^{\gamma-1} \, \nabla u_\gamma) + \frac{1}{\gamma^\alpha} u_\gamma  \, \Phi(p_\gamma).
\end{aligned}
\end{equation}
Hence we multiply equation \eqref{EQUGAMMA} by $\varphi^\prime_\delta(u_\gamma)$, where for each $\delta>0$, 
$$
\varphi_\delta(z)= 
  \left \{
\begin{aligned}
  & (z^2 + \delta^2)^{1/2} - \delta & \quad \text{for $z \geq 0$}, 
  \\
  & 0 & \quad  \text{for $z \leq 0$},
\end{aligned}
\right.
$$
and using integration by parts we obtain 
$$
\begin{aligned}
\int_{\Omega} \partial_t \varphi_\delta(u_\gamma) \ dx &
    = \gamma \int_{\Omega} \text{div}(\Theta(u_\gamma) \, u_\gamma^{\gamma-1} \nabla u_\gamma) \, \varphi^\prime_\delta(u_\gamma) \ dx 
   + \frac{1}{\gamma^\alpha} \int_{\Omega} u_\gamma \, \Phi(p_\gamma) \, \varphi^\prime_\delta(u_\gamma) \ dx 
   \\[5pt]
   & \leq - \int_{\Omega} \Theta(u_\gamma) \, u_\gamma^{\gamma-1} |\nabla u_\gamma|^2 \varphi^{\prime \prime}_\delta(u_\gamma) \, dx + \int_{\Gamma} \Theta(u_\gamma) \varphi^\prime_\delta(u_\gamma)  \nabla p_{\gamma}  \cdot \mathbf{n} \ dr 
   \\[5pt]
   & + \frac{\Phi(0)}{\gamma^\alpha} \int_{\Omega} u_\gamma \, \varphi^\prime_\delta(u_\gamma) \ dx. 
\end{aligned}
$$
Letting $\delta \to 0^+$ and provided that $\Theta(\cdot) \geq 0$, we arrive to
\begin{align} 
\label{eq: L1 contraction}
\frac{d}{dt} \int_{\Omega} |u_\gamma|^{+} \ dx \leq  \int_{\Gamma} |f_\gamma| \ dr + \frac{\Phi(0)}{\gamma^\alpha} \int_{\Omega} |u_\gamma|^+(u_\gamma) \ dx, 
\end{align}
where we have used that $\varphi^{\prime\prime}_\delta(\cdot)> 0$. Gronwall's inequality yields for any $t \in [0,T]$ that
\begin{align*}
\int_{\Omega} |u_\gamma|^{+} \ dx \leq e^{\frac{\Phi(0)T}{\gamma^\alpha}} \big( \int_{\Omega} |u_{0\gamma}|^{+} \ dx
+ \int_{\Gamma_T} |f_\gamma| \ drdt \big). 
\end{align*}
Analogously, we derive a similar estimate for the $L^1$-norm of $|u|^{-}$. This ends the proof of (\ref{eq: L1 estimate ugamma}).

\medskip
2. Now, for every $\gamma>1$, we write $p_\gamma(u_\gamma)=  u_\gamma \, u_\gamma^{\gamma-1}$.  
Then, the comparison given in (\ref{eq: comparison}) yields $0 \leq u_\gamma \leq (p_M)^{1/\gamma}$ a.e. in $Q_T$.  
Hence, we obtain (\ref{eq: L1 estimate pgamma}) for every $t \in [0,T]$.
\end{proof}

\subsection{Benilan-Aronson type estimates.}
\label{Benilan-Aronson type estimates}

We draw upon the idea presented in \cite{aronson1979} to establish the proof for the following  

\begin{proposition} \label{prop: Benilan-aronson estimates}
Let us assume that the set of Hypotheses \ref{condition: coefficients of the equations}-\ref{condition: initial condition for the compressible BL} hold.  Fix $\alpha>1$ 
and for every $\gamma>1$ let $u_\gamma$ be the unique solution of \eqref{eq: the compressible BL system} in the sense of Theorem \ref{thm: existence uniq strong solution of the compressible BL}.  
Then, for every $\gamma>1$ and a.e. in $Q_T$, it follows that
\begin{align} \label{eq: Benilan-Aronson}
\Theta(u_\gamma) \, \Delta p_\gamma + \frac{1}{\gamma^\alpha} u_\gamma \, \Phi(p_\gamma) \geq - \frac{r_\Phi}{\gamma^\alpha} u_\gamma \, \frac{e^{- \frac{r_\Phi}{\gamma^{\alpha-1}}t}}{1- e^{- \frac{r_\Phi}{\gamma^{\alpha-1}}t}}, 
\end{align}
where $r_\Phi:= \displaystyle \min \Big \{  \Phi(p) - \Phi'(p) p \Big \}$, which is a non-negative function by hypothesis done on $\Phi$.
\end{proposition}

\begin{proof}
1. First, it is enough to consider $u_\gamma> 0$. Indeed, if $u_\gamma= 0$, then the inequality \eqref{eq: Benilan-Aronson} is trivially satisfied. 
Therefore, let us consider the equation for the pressure $p_\gamma$, that is, we multiply equation \eqref {EQUGAMMA} by $\gamma \, u_\gamma^{\gamma-1}$ to obtain 
 \begin{equation} 
\label{eq: BA1}
\begin{aligned}
   \partial_t p_\gamma &= \gamma \, \text{div} \Big ( \Theta(u_\gamma) \nabla p_\gamma \Big ) u_\gamma^{\gamma-1}+ \frac{1}{\gamma^{\alpha-1}} p_\gamma \Phi(p_\gamma) 
   \\
   & =  \gamma \Big (\nabla \Theta(u_\gamma) \cdot  \nabla p_\gamma + \Theta(u_\gamma) \Delta p_\gamma \Big ) u_\gamma^{\gamma-1}+ \frac{1}{\gamma^{\alpha-1}} p_\gamma \Phi(p_\gamma). 
\end{aligned}
\end{equation} 
Here, we observe that
\begin{align*}
\gamma u_\gamma^{\gamma-1} (\nabla \Theta(u_\gamma) \cdot \nabla p_\gamma)= \gamma \Theta^\prime(u_\gamma) u^{\gamma-1}_\gamma \nabla u_\gamma \cdot \nabla p_\gamma = \Theta^\prime(u_\gamma) |\nabla p_\gamma|^2.
\end{align*}
Hence the equation \eqref{eq: BA1} reads as
\begin{align*}
\partial_t p_\gamma - \Theta^\prime(u_\gamma) |\nabla p_\gamma|^2 
- \gamma \, \Theta(u_\gamma) \Delta p_\gamma \frac{u_\gamma^{\gamma}}{u_\gamma}= \frac{1}{\gamma^{\alpha-1}} p_\gamma \Phi(p_\gamma),
\end{align*}
which is equivalent to
\begin{equation}
\label{EQBNL}
\partial_t p_\gamma = \Theta^\prime(u_\gamma) |\nabla p_\gamma|^2 + \gamma p_\gamma \Big ( \frac{\Theta(u_\gamma)}{u_\gamma} \Delta p_\gamma + \frac{1}{\gamma^\alpha}\Phi(p_\gamma) \Big ).
\end{equation}

\medskip
2. To follow, we write 
\begin{equation}
\label{eq: identity w}
w=  \frac{\Theta(u_\gamma)}{u_\gamma} \Delta p_\gamma +\frac{1}{\gamma^\alpha} \Phi(p_\gamma) =: \bar{\Theta}(u_\gamma) \Delta p_\gamma + \frac{1}{\gamma^{\alpha}} \Phi(p_\gamma), 
\end{equation}
and using the summation convention, we have from equation \eqref{EQBNL}, 
\begin{align*}
\partial_t p_\gamma =  \Theta^\prime(u_\gamma) \frac{\partial p_\gamma}{\partial x_k} \frac{\partial p_\gamma}{\partial x_k} + \gamma \, p_\gamma \, w.
\end{align*}
Then, we differentiate in order to $x_i$ , $(i=1, \ldots, d)$, the above equation to obtain
\begin{align*}
\partial_t \frac{\partial p_\gamma}{\partial x_i}
= \frac{\partial  \Theta^\prime(u_\gamma)}{\partial x_i} \frac{\partial p_\gamma}{\partial x_k} \frac{\partial p_\gamma}{\partial x_k} 
+ 2  \Theta^\prime(u_\gamma) \frac{\partial p_\gamma}{\partial x_k} \frac{\partial^2 p_\gamma}{\partial x_k \partial x_i} 
+ \gamma \frac{\partial p_\gamma}{\partial x_i} w +\gamma p_\gamma \frac{\partial w}{\partial x_i},
\end{align*}
and again differentiating in order to $x_j$, we have
\begin{align*}
\partial_t \frac{\partial^2 p_\gamma}{\partial x_i \partial x_j}
&= \frac{\partial^2  \Theta^\prime(u_\gamma)}{\partial x_i \partial x_j} |\nabla p_\gamma|^2 
+ 2 \frac{\partial  \Theta^\prime(u_\gamma)}{\partial x_i} \frac{\partial^2 p_\gamma}{\partial x_k \partial x_j} \frac{\partial p_\gamma}{\partial x_k} 
+ 2 \frac{\partial  \Theta^\prime(u_\gamma)}{\partial x_j} \frac{\partial p_\gamma}{\partial x_k} \frac{\partial^2 p_\gamma}{\partial x_k \partial x_i} 
\\
&+ 2 \Theta^\prime(u_\gamma) \frac{\partial^2 p_\gamma}{\partial x_k \partial x_j} \frac{\partial^2 p_\gamma}{\partial x_k \partial x_i} 
+ 2  \Theta^\prime(u_\gamma) \frac{\partial p_\gamma}{\partial x_k} \frac{\partial^3 p_\gamma}{\partial x_k \partial x_i \partial x_j} 
+\gamma \frac{\partial^2 p_\gamma}{\partial x_i \partial x_j} w 
\\
& + \gamma \frac{\partial p_\gamma}{\partial x_i} \frac{\partial w }{\partial x_j} 
 + \gamma \frac{\partial p_\gamma}{\partial x_j} \frac{\partial w}{\partial x_i} 
+ \gamma  p_\gamma \frac{\partial^2 w}{\partial x_i \partial x_j}.
\end{align*}
\no Contracting $i=j$ implies
$$
\begin{aligned}
\partial_t \Delta p_\gamma &= \Delta \Theta^\prime(u_\gamma) |\nabla p_\gamma|^2 
+ 4 \nabla \Theta^\prime(u_\gamma) \cdot (D^2 p_\gamma) \nabla p_\gamma 
\\
&+ 2  \Theta^\prime(u_\gamma) (D^2 p_\gamma) : (D^2 p_\gamma)
+ 2  \Theta^\prime(u_\gamma) \nabla (\Delta p_\gamma) \cdot \nabla p_\gamma 
\\
&+\gamma \Delta p_\gamma w 
+ 2 \gamma \nabla p_\gamma \cdot \nabla w 
+ \gamma p_\gamma \Delta w,
\end{aligned}
$$
and setting $v= \Delta p_\gamma$, we obtain
\begin{align}  \label{eq: BA2}
\partial_t v &= \Delta  \Theta^\prime(u_\gamma) |\nabla p_\gamma|^2 + 4 \nabla  \Theta^\prime(u_\gamma) \cdot (D^2 p_\gamma) \nabla p_\gamma
+ 2 \Theta^\prime(u_\gamma) (D^2 p_\gamma) : (D^2 p_\gamma)  \nonumber 
\\
&+ 2 \Theta^\prime(u_\gamma) \nabla v \cdot  \nabla p_\gamma 
 + \gamma v w + 2 \gamma  \nabla p_\gamma \cdot \nabla w +  \gamma p_\gamma \Delta w.
\end{align}

\medskip
3. Now, we recall the equation for the pressure term $p_\gamma$, that is, 
\begin{align*}
\partial_t p_\gamma = \Theta^\prime(u_\gamma) |\nabla p_\gamma|^2
 + \gamma p_\gamma w,
 \end{align*}
and differentiate equation \eqref{eq: identity w} with respect to time, we obtain
\begin{align} \label{eq: equation for w}
\partial_t w &=  \partial_t \bar \Theta(u_\gamma) v +  \bar \Theta(u_\gamma) \partial_t v + \frac{1}{\gamma^\alpha} \Phi'(p_\gamma) \partial_t p_\gamma.
\end{align}
Moreover, we assume without loss of generality that $\Theta^\prime(u_\gamma)> 0$. Otherwise, we replace $\Theta$ by 
$$
  \Theta_\gamma(z)=  \Theta(z) + \frac{1}{\gamma^\alpha} z.
$$  
Hence we obtain from equation \eqref{eq: BA2} that
\begin{align} 
\label{eq: BA_ineq_v}
\partial_t v &\geq \Delta \Theta^\prime(u_\gamma) |\nabla p_\gamma|^2 + 4 \nabla \Theta^\prime(u_\gamma) \cdot (D^2 p_\gamma) \nabla p_\gamma \nonumber \\
& + 2 \Theta^\prime(u_\gamma) \nabla v \cdot  \nabla p_\gamma + \gamma v w + 2 \gamma \nabla p_\gamma \cdot \nabla w \nonumber \\
& +\gamma  p_\gamma \Delta w.
\end{align}
In what follows we calculate
$\partial_t \Phi(p_\gamma) = \Phi'(p_\gamma) \partial_t p_\gamma$, 
and use \eqref{eq: equation for w}, \eqref{eq: BA_ineq_v}, that is to say
\begin{equation} 
\label{eq: BA_final}
\begin{aligned} 
\partial_t w & \geq \partial_t \bar \Theta(u_\gamma) v 
+ \bar \Theta(u_\gamma) \Big ( \Delta \Theta^\prime(u_\gamma) |\nabla p_\gamma|^2 
+ 4 \nabla \Theta^\prime(u_\gamma) \cdot (D^2 p_\gamma) \nabla p_\gamma
\\
&+ 2 \Theta^\prime(u_\gamma) \nabla v \cdot \nabla p_\gamma 
+ 2 \nabla p_\gamma \cdot \nabla w \Big )  
+ \gamma \bar \Theta(u_\gamma) v w +\gamma  \bar \Theta(u_\gamma)  p_\gamma \Delta w  
\\
& + \frac{1}{\gamma^\alpha} \Phi'(p_\gamma) \Theta^\prime(u_\gamma) |\nabla p_\gamma|^2 + \Phi'(p_\gamma) p_\gamma w. 
\end{aligned}
\end{equation}
From \eqref{eq: identity w} it is immediate that
\begin{align*}
\gamma \, \bar \Theta(u_\gamma) \, v \, w= \gamma w^2 - \frac{1}{\gamma^{\alpha-1}} \Phi(p_\gamma) w, 
\end{align*}
and defining  
\begin{align*}
F_\gamma(u_\gamma, p_\gamma):= 4 \nabla \Theta^\prime(u_\gamma) \cdot (D^2 p_\gamma) \nabla p_\gamma 
+ 2 \Theta^\prime(u_\gamma) \nabla v \cdot \nabla p_\gamma + 2 \gamma \nabla p_\gamma \cdot \nabla w, 
\end{align*}
we have from \eqref{eq: BA_final} that 
\begin{align*}
\partial_t w & \geq \Theta^\prime(u_\gamma) F_\gamma(u_\gamma, p_\gamma) + \partial_t \Theta^\prime(u_\gamma) v +\gamma  \Theta^\prime(u_\gamma) v w 
\\
& + \Theta^\prime(u_\gamma) \gamma p_\gamma w + \frac{1}{\gamma^\alpha} \Theta^\prime(p_\gamma) |\nabla p_\gamma|^2 \Theta^\prime(u_\gamma)
+ \frac{1}{\gamma^\alpha} \Phi'(p_\gamma) \gamma p w.
\end{align*}
Therefore, noting that $p_\gamma \geq 0$ and $\Theta^\prime(u_\gamma)>0$, it follows that
\begin{equation}
\label{INEQ}
\begin{aligned}
\partial_t w &\geq \Theta^\prime(u_\gamma) F_\gamma(u_\gamma, p_\gamma) + \partial_t \Theta^\prime(u_\gamma) v 
+ \gamma w^2 - \frac{w}{\gamma^{\alpha-1}} (\Phi(p_\gamma) - \Phi'(p_\gamma) p_\gamma) 
\\
&+ \frac{\Theta^\prime(u_\gamma)}{\gamma^\alpha} \nabla \Phi(p_\gamma) \cdot \nabla p_\gamma.
\end{aligned}
\end{equation}

\medskip
4. Finally, let us search for a particular solution of the equation \eqref{INEQ}, with $p_\gamma(t,x)= P_\gamma(t)$, that is independent of $x$. 
Therefore, $u_\gamma= U_\gamma(t)$ for certain function $U_\gamma$ independent of $x$, we have $F_\gamma(U_\gamma,P_\gamma)=0$, and $v=V_\gamma(t)=0$. 
Then, we consider the following Ricatti equation for $w= W(t)$, that is to say, 
\begin{align*}
\dot W(t)= \gamma W^2(t) - \frac{1}{\gamma^{\alpha-1}} r_\Phi W(t), 
\end{align*}
where $r_\Phi:= \displaystyle \min_{0 \leq p \leq p_M} (\Phi(p)- \Phi'(p)p) >0$ by hypothesis.
Setting $Z(t)=W^{-1}(t)$, we derive the following equation
\begin{align*}
\dot Z(t)= - \frac{\gamma W^2(t)}{W^2(t)} + \frac{r_\Phi W(t)}{\gamma^{\alpha-1} W(t)}
\end{align*}
or equivalently
\begin{align*}
\dot Z(t) -  \frac{r_\Phi}{\gamma^{\alpha-1}} Z(t)= - \gamma.
\end{align*}
\no Multiplying the equation above by $e^{- \frac{r_\Phi}{\gamma^{\alpha-1}}t}$ yields
\begin{align*}
\frac{d}{dt}\Big ( Z(t) e^{- \frac{r_\Phi}{\gamma^{\alpha-1}t}} \Big )= - \gamma e^{- \frac{r_\Phi}{\gamma^{\alpha-1}}t}.
\end{align*}
\no Integrating on $[0,t]$ and assuming $Z(0)=0$, which implies that $W$ explodes at $0$, it follows that 
\begin{align*}
Z(t)e^{- \frac{r_\Phi}{\gamma^{\alpha-1}}t}= - \frac{\gamma^\alpha}{r_\Phi}(1- e^{- \frac{r_\Phi}{\gamma^{\alpha-1}}t}),
\end{align*}
and hence
\begin{align*}
W(t) = - \frac{r_\Phi}{\gamma^\alpha} \frac{e^{- \frac{r_\Phi}{\gamma^{\alpha-1}}t}}{1- e^{- \frac{r_\Phi}{\gamma^{\alpha-1}}t}}.
\end{align*}
Therefore, we have 
\begin{align*}
w(t)= \frac{\Theta(u_\gamma)}{u_\gamma} \Delta p_\gamma + \frac{1}{\gamma^\alpha} \Phi(p_\gamma) \geq
 - \frac{r_\Phi}{\gamma^\alpha} \frac{e^{- \frac{r_\Phi}{\gamma^{\alpha-1}}t}}{1- e^{- \frac{r_\Phi}{\gamma^{\alpha-1}}t}},
\end{align*}
from which follows the desired Benilan-Aronson estimate. 
\end{proof}

\begin{corollary}
 \label{corollary: pointwise positive limit time derivatives} 
 Under the assumptions of Theorem \ref{thm: main theorem}, one has that $\partial_t u_\infty, \partial_t p_\infty \geq 0$
in the sense of distributions. 
\end{corollary}
\begin{proof}
\no Recalling the equation for the pressure $p_\gamma$, that is,  
\begin{align*}
\partial_t p_\gamma = \Theta^\prime(u_\gamma) |\nabla p_\gamma|^2 + \gamma p_\gamma w,
\end{align*}
and since $p_\gamma \in [0,p_M]$, $\Theta^\prime(u_\gamma) >0$, we have 
\begin{align*}
\partial_t p_\gamma & \geq - \frac{r_\Phi}{\gamma^{\alpha-1}} \frac{e^{- \frac{r_\Phi t}{\gamma^{\alpha-1}}}}{1- e^{- \frac{r_\Phi t }{\gamma^{\alpha-1}}}} p_\gamma. 
\end{align*}
Then, passing to the limit as $\gamma \to \infty$ in the above estimate, we obtain 
$\partial_t p_\infty \geq 0$ in distribution sense, where we have used Proposition \ref{prop_ L1 bounds pgamma, ugamma}. 

Now, we observe that 
\begin{equation} 
 \label{eq: pointwise time derivative for u}
 \begin{aligned}
\partial_t u_\gamma &= \partial_t p_\gamma^{\frac{1}{\gamma}} = \frac{1}{\gamma} p_\gamma^{\frac{1}{\gamma}-1} \partial_t p_\gamma=  \frac{1}{\gamma}  \frac{1}{p_\gamma^{1- \frac{1}{\gamma}}} \partial_t p_\gamma
\\[5pt]
& \geq \frac{1}{p_M^{1- \frac{1}{\gamma}}} \times \frac{(- r_\Phi)}{\gamma^{\alpha}} \frac{e^{- \frac{r_\Phi t }{\gamma^{\alpha-1}}}}{1 - e^{- \frac{r_\Phi t}{\gamma^{\alpha-1}}}} p_\gamma, 
\end{aligned}
\end{equation}
from which, similarly, we conclude that $\partial_t u_\infty \geq 0$ in distribution sense. 
\end{proof}

\subsection{Uniform bounds for the time derivatives.}

\begin{proposition}
\label{eq: prop uniform bound time derivatives} 
Under the asumptions of Theorem \ref{thm: main theorem}, there exists a positive constant $C= C(T)$ such that, for $\gamma> 1$ with $1/\gamma< T$, 
\begin{align*}
||\partial_t u_\gamma||_{L^1(Q_T^\gamma)} + || \partial_t p_\gamma||_{L^1(Q_T^\gamma)} \leq C,
\end{align*}
where $Q_T^\gamma= (\gamma^{-1},T) \times \Omega$. 
\end{proposition}
\begin{proof}
1. First, using the fact that $|z|= z + 2 |z|^{-}$, equations \eqref{eq: L1 contraction}, \eqref{eq: pointwise time derivative for u}, and $0 \leq u_\gamma \leq 1$, we may write for every $t \in [0,T]$ 
\begin{align*}
\int_{ \Omega} |\partial_t u_\gamma(t)| \, dx  &= \partial_t \int_{ \Omega} u_\gamma(t) \, dx + 2 \int_{ \Omega} |\partial_t u_\gamma(t)|^- \, dx 
\\
& \leq \int_{\Gamma} |f_\gamma(t)| \ dr + \frac{\Phi(0)}{\gamma^{\alpha}} \int_{ \Omega} |u_\gamma(t)| \, dx 
+ 2 \frac{1}{p_M^{1- \frac{1}{\gamma}}} \frac{r_\Phi}{\gamma^{\alpha}}  \int_{ \Omega} \frac{e^{- \frac{r_\Phi t }{\gamma^{\alpha-1}}}}{1 - e^{- \frac{r_\Phi t}{\gamma^{\alpha-1}}}} \, p_\gamma(t) \, dx. 
\end{align*}
Due to (\ref{eq: L1 estimate ugamma}), $0 \leq u_\gamma(t) \leq 1$ and the above estimate, it follows that
\begin{equation} 
\label{eq: estimate final for derivative in time}
\begin{aligned}
\int_{\frac{1}{\gamma}}^T \!\!\! \int_{ \Omega} |\partial_t u_\gamma| dx dt 
& \leq  C(T)
+ \frac{2 r_\Phi}{P_M^{1- \frac{1}{\gamma}} \gamma^\alpha} \int_{\frac{1}{\gamma}}^T \Big (  \frac{e^{- \frac{r_\Phi t}{\gamma^{\alpha-1}}}}{1- e^{- \frac{r_\Phi t }{\gamma^{\alpha-1}}}} \int_{\Omega} u_\gamma u_\gamma^{\gamma-1} dx  \Big )dt 
\\[5pt]
& \leq C(T) + \frac{2 r_\Phi}{P_M^{1- \frac{1}{\gamma}} \gamma^\alpha} e^{\frac{\Phi(0)T}{\gamma^\alpha}} \int_{\frac{1}{\gamma}}^T \frac{e^{- \frac{r_\Phi t}{\gamma^{\alpha-1}}}}{1- e^{- \frac{r_\Phi t}{\gamma^{\alpha-1}}}} dt
\\[5pt]
& \leq C(T) +  \frac{2 r_\Phi}{P_M^{1- \frac{1}{\gamma}} \gamma^\alpha} e^{\frac{\Phi(0)T}{\gamma^\alpha}} 
\frac{\gamma^{\alpha-1}}{r_\Phi}  \big ( \ln \big ( 1 - e^{- \frac{r_\Phi }{\gamma^{\alpha-1}}T} \big ) - \ln \big ( 1- e^{- \frac{r_\Phi }{\gamma^{\alpha-1}}\frac{1}{\gamma}} \big ) \big )
\\[5pt]
& \leq C(T) + C_1(T) \frac{1}{\gamma} \big ( \ln \big ( 1 - e^{- \frac{r_\Phi }{\gamma^{\alpha-1}}T} \big ) - \ln \big ( 1- e^{- \frac{r_\Phi }{\gamma^{\alpha-1}}\frac{1}{\gamma}} \big ) \big )  
\\[5pt]
& \leq C_2(T)< \infty. 
\end{aligned}
\end{equation} 

\medskip
2. Now, since $p_\gamma(u_\gamma)= u_\gamma^\gamma$ we have,
$|\partial_t p_\gamma| = |\partial_t u_\gamma^\gamma| = \gamma \, u_\gamma^{\gamma-1} \, |\partial_t u_\gamma|$. 
Therefore, we obtain 

\begin{align} \label{eq: lim pgamma derivative}
\int_{\frac{1}{\gamma}}^T \!\!\! \int_{\Omega} |\partial_t p_\gamma| dx dt 
& \leq C(T) + C_1(T) \big ( \ln \big ( 1 - e^{- \frac{r_\Phi }{\gamma^{\alpha-1}}T} \big ) - \ln \big ( 1- e^{- \frac{r_\Phi }{\gamma^{\alpha-1}}\frac{1}{\gamma}} \big ) \big )  
\\
& \leq C_2(T)< \infty,
\end{align}

where we have used that, $0 \leq u_\gamma(t,x) \leq 1$ and equation \eqref{eq: estimate final for derivative in time}, which concludes the proof. 
\end{proof}

\no Taking the supremum in $\gamma >1$ in (\ref{eq: estimate final for derivative in time}) and (\ref{eq: lim pgamma derivative}) yields the following result.

\begin{proposition} \label{prop: prop uniform bound time derivatives conclusion} 
Under the asumptions of Theorem \ref{thm: main theorem}, there exists a positive constant $C= C(T)$ such that, for $\gamma> 1$ with $1/\gamma< T$, 
\begin{align*}
||\partial_t u_\gamma||_{L^1(Q_T)} + || \partial_t p_\gamma||_{L^1(Q_T)} \leq C.
\end{align*}
\end{proposition}

\subsection{A priori estimates for the spatial gradients.}
\begin{proposition} \label{prop: apriori estimates gradients}
Under the assumptions of Theorem \ref{thm: main theorem}, there exists $C=C(T)> 0$ such that, for every $\gamma>1$ it holds
\begin{align*}
||\nabla u_\gamma||_{L^1(Q_T)} + ||\nabla  p_\gamma||_{L^1(Q_T)} \leq C.
\end{align*}
\end{proposition}
\begin{proof}
1. First, let us estimate $||\nabla u_\gamma||_{L^1(Q_T)}$. For that 
we recall the equation satisfyied  for $u_\gamma$ given by $(\ref{eq: the system for comparison})_1$, that is to say
\begin{align*}
\partial_t u_\gamma = \Delta \Psi(u_\gamma) + \frac{1}{\gamma^\alpha} u_\gamma \Phi(p_\gamma),
\end{align*}
and let us denote $w_{i \gamma}= \partial_{x_i} u_\gamma$ , $(i= 1, \ldots, d)$, 
$\mathbf{w}_\gamma \equiv (w_{1 \gamma}, \ldots, w_{d \gamma})$.
Then we obtain from the above equation 
\begin{equation}
\label{EQGRADEST}
\partial_t w_{i \gamma} = \Delta(\tilde{w}_{i\gamma}) 
+ \frac{\Phi(p_\gamma)}{\gamma^{\alpha}} \, w_{i \gamma} 
+ \frac{1}{\gamma^{\alpha}} u_\gamma \,  \Phi'(p_\gamma) \partial_{x_i} p_\gamma, \quad (i= 1, \ldots, d), 
\end{equation}
where $\tilde{w}_{i\gamma}= \Psi^\prime(u_\gamma) w_{i\gamma}$.  

\medskip
2. Now, we define $$-\beta:= \displaystyle \min_{0 \leq p \leq p_M} \Phi'(p)< 0,$$
and multiply the $i$th equation of \eqref{EQGRADEST} by 
$\frac{\partial}{\partial z_i} \vp_\delta(\mathbf{w}_\gamma)$, where $\vp_\delta(\mathbf{z})= (|\mathbf{z}|^2 + \delta^2)^{1/2}$,
hence we add it up and integrate on $(0,t) \times \Omega$, for any $0< t \leq T$. 
Then, considering the main resulting terms
$$
\begin{aligned}
\quad   \int_0^t \!\!\! \int_\Omega  \frac{\partial w_{i\gamma}}{\partial t} \, \frac{\partial \vp_\delta(\mathbf{w}_\gamma)}{\partial z_i}  \, dx d\tau
    = \int_\Omega \vp_\delta(\mathbf{w}_\gamma(t)) \, dx - \int_\Omega \vp_\delta(\mathbf{w}_\gamma(0)) \, dx, \hspace{1000pt} 
\end{aligned}    
$$
and 
$$
\begin{aligned}
\quad   \int_0^t \!\!\! \int_\Omega \frac{\partial}{\partial x_j} &  \Big(\frac{\partial \tilde{w}_{i\gamma}}{\partial x_j}\Big) \, \frac{\partial \vp_{\tilde{\delta}}(\tilde{\mathbf{w}}_\gamma)}{\partial z_i}  \, dx d\tau
\\[5pt]
& = \int_0^t \!\!\! \int_\Gamma \ \frac{\partial \tilde{w}_{i\gamma}}{\partial x_j} \, n_j \, \frac{\partial \vp_{\tilde{\delta}}(\tilde{\mathbf{w}}_\gamma)}{\partial z_i} \, dr d\tau 
-  \int_0^t \!\!\! \int_\Omega \frac{\partial \tilde{w}_{i\gamma}}{\partial x_j} \, \frac{\partial^2 \vp_{\tilde{\delta}}(\tilde{\mathbf{w}}_\gamma)}{\partial z_i \partial z_k} 
\, \frac{\partial \tilde{w}_{k\gamma}}{\partial x_j} \, dx d\tau, \hspace{1000pt}
\end{aligned}
$$
we obtain passing to the limit as $\delta \to 0^+$
\begin{equation}
\label{GRADINEQ} 
\begin{aligned}
\int_\Omega |\mathbf{w}_\gamma(t)| \, dx & +  \beta \, \int_0^t \!\!\! \int_\Omega u_\gamma \, |\nabla p_\gamma| \, dxd\tau 
\\[5pt]
& \leq \int_\Omega |\mathbf{w}_\gamma(0)| \, dx + \frac{\Phi(0)}{\gamma^{\alpha}}  \int_0^t \!\!\! \int_\Omega  |\mathbf{w}_\gamma| \, dxd\tau 
+ \int_0^t \!\!\! \int_\Gamma \ \big| \nabla |\nabla \Psi(u_\gamma)| \cdot \mathbf{n}\big| \, dr d\tau, 
\end{aligned}
\end{equation}
where we have used that $0< \Psi^\prime(u_\gamma) \leq C$, for some positive constant $C$. 

\medskip
3. To follow, since $\Gamma$ has $C^2$ regularity, $\mathbf{n}$ is a $C^1$ vector value function and also could be extended
to a function $\mathbf{N}$, defined on the whole of $\RR^d$ such that, $\mathbf{N} \in C^1(\RR^d; \RR^d)$ with $|\mathbf{N}| \leq 1$.
Therefore, we may write 
$$
    \nabla \Psi(u_\gamma) \cdot \mathbf{N}= f_\gamma \quad \text{on $\Gamma_T$}, 
$$
where $f_\gamma(t,\cdot)$ is also defined in $\RR^d$ for each $t \in (0,T)$.  Moreover, we have on $\Gamma_T$ that 
\begin{equation}
\label{ESTIMBORDO1}
  \big| \nabla |\nabla \Psi(u_\gamma)| \cdot \mathbf{N}\big| \leq \frac{\sqrt{d}}{|\nabla \Psi(u_\gamma) \cdot \mathbf{N}|}
  | D^2  \Psi(u_\gamma) \, \mathbf{N} | 
  \leq \frac{\sqrt{d}}{f_b} \
  | D^2  \Psi(u_\gamma) \, \mathbf{N} |. 
\end{equation}
Due to $\nabla \Psi(u_\gamma) \cdot \mathbf{N} - f_\gamma= 0$ on $\Gamma_T$, we have that 
$$
  \big| \nabla \big(\nabla \Psi(u_\gamma) \cdot \mathbf{N} - f_\gamma \big)\big|
  = \big|\nabla \big(\nabla \Psi(u_\gamma) \cdot \mathbf{N} - f_\gamma \big) \cdot \mathbf{N} \big|
   \quad \text{on $\Gamma_T$},
$$
and thus it follows on $\Gamma_T$ that 
\begin{equation}
\label{ESTIMBORDO2}
  | D^2  \Psi(u_\gamma) \, \mathbf{N} | \leq | \nabla f_\gamma| + K \, |f_\gamma|, 
\end{equation}
where $K(r)= \sum_{j=1}^{d-1} |\kappa_j(r)|$, and $\kappa_1, \ldots, \kappa_{d-1}$ are the principal curvatures of $\Gamma$ at $r$, 
(see Gilbarg, Trudinger \cite{gilbarg1977elliptic}, p. 354). 
Then, from \eqref{GRADINEQ}, \eqref{ESTIMBORDO1}, and \eqref{ESTIMBORDO2} we obtain
$$
  \begin{aligned}
\int_\Omega |\nabla u_\gamma(t)| \, dx  
 \leq \int_\Omega |\nabla u_{\gamma,0}| \, dx 
 &+ \frac{\sqrt{d}}{f_b} \int_0^t \!\!\! \int_\Gamma \big(| \nabla f_\gamma| + K \, |f_\gamma| \big)\, dr d\tau
 \\[5pt]
 &+ \frac{\Phi(0)}{\gamma^{\alpha}}  \int_0^t \!\!\! \int_\Omega  |\nabla u_\gamma| \, dxd\tau.
\end{aligned}
$$
Gronwall's inequality allow us to conclude that
\begin{align*}
\iint_{Q_T} |\nabla u_\gamma| \, dxdt 
\leq C_1(T) \quad \text{and also} \quad \beta \iint_{Q_T} u_\gamma \, |\nabla p_\gamma| \, dxdt \leq C_1(T), 
\end{align*}
where $C_1(T)$ is a positive constant, and  we have used \eqref{BONDGRADU0}, \eqref{BONDFLUX}. 

\medskip
4. Finally, we estimate $||\nabla p_\gamma||_{L^1(Q_T)}$.
From \eqref{stiffpressure}, we have $\nabla p_\gamma= \gamma \, u_\gamma^{\gamma-1} \nabla u_\gamma$. 
For $\gamma> 1$, it follows that $\gamma \, u_\gamma^{\gamma-1} \leq 1$ on $\{u_\gamma \leq 1/2\}$, 
and similarly $1 \leq 2 \, u_\gamma$ on $\{u_\gamma \geq 1/2\}$. Therefore, we may write 
$$
\begin{aligned}
 \iint_{Q_T} |\nabla p_\gamma| \, dxdt &\leq \iint_{Q_T \cap \{u_\gamma \leq 1/2\}} \gamma \, u_\gamma^{\gamma-1} |\nabla u_\gamma| \, dxdt 
 + 2 \iint_{Q_T \cap \{u_\gamma \geq 1/2\}} u_\gamma \, |\nabla p_\gamma| \, dxdt 
 \\[5pt]
 &\leq C_2(T), 
\end{aligned}
$$
where $C_2(T)$ is a positive constant, and thus the proof is finished. 
\end{proof}

\subsection{Compactness argument and properties of the limit.}

\begin{theorem} 
\label{them: theorem convergence} 
Under the assumptions of Theorem \ref{thm: main theorem} we have that the 
family of solutions $\{(u_\gamma, p_\gamma)\}$ of 
the system \eqref{stiffpressure}, \eqref{eq: the compressible BL system} and \eqref{BCIC}
converges strongly in $L^1(Q_T)$ 
as $\gamma \ra \infty$ to an element $(u_\infty, p_\infty)$. Moreover, one has that, $u_\infty,  p_\infty \in BV(Q_T)$.
\end{theorem}

\begin{proof}
1. First, we observe that the previous sections work culminates on the conclusion that, $\{u_\gamma\}_{\gamma>1}$ 
and $\{p_\gamma\}_{\gamma>1}$ are uniformly (with respect to $\gamma>0$) bounded in $W^{1,1}(Q_T)$.  
Therefore,  passing to a subsequence (we maintain the same notation),
we have strong convergence in $L^1(Q_T)$ for both $\{u_\gamma\}$ and $\{p_\gamma\}$ since $Q_T$ is bounded. 


\medskip
2. Finally, the fact that $u_\infty,  p_\infty \in BV(Q_T)$ follows from standard arguments 
of weak convergence and uniqueness of the limit.
\end{proof}

\medskip
The following proposition tell us in particular that $u_\infty$, $p_\infty$ satisfy 
the Hele-Shaw limit.
\begin{proposition} 
\label{prop: Hele Shaw limit}
Under the assumptions of Theorem \ref{them: theorem convergence}, one has that \eqref{eq: comparison} yields
\begin{align}
 \label{eq: Hele-Shaw limit}
0 \leq u_\infty \leq 1 \quad \text{ and } \quad 0 \leq p_\infty \leq p_M \quad \text{a.e.}-Q_T.
\end{align}
Moreover the functions $u_\infty$ and $p_\infty$ satisfy the Hele-Shaw graph property
 $p_\infty (u_\infty -1)= 0$.
\end{proposition}

\begin{proof}
Let $\{(u_\gamma, p_\gamma)\}$ be the family of solutions of 
the system \eqref{stiffpressure}, \eqref{eq: the compressible BL system} and \eqref{BCIC}, which
converges strongly in $L^1(Q_T)$ due to Theorem \ref{them: theorem convergence}. Then, 
after extraction of subsequences, we can pass to the almost everywhere limit in the equation \eqref{eq: comparison},
that is,
$$
0 \leq u_\gamma \leq \sqrt[\gamma]{p_M} \quad \text{ and } \quad 0 \leq p_\gamma \leq p_M \quad \text{a.e.}-Q_T, 
$$
to obtain \eqref{eq: Hele-Shaw limit}. 
Similarly, we can pass to the almost everywhere limit in the equation 
$$
u_\gamma \, p_\gamma= \big(p_\gamma\big)^{(\gamma+1)/\gamma}
$$
to obtain that $p_\infty (u_\infty -1)= 0$. 
\end{proof}

\subsection{$L^2$ bounds for $ \nabla p_\gamma$ and weak convergence.}

\begin{theorem} 
\label{thm: l2 bounds on the gradient of the pressure}
Under the assumptions of Theorem \ref{them: theorem convergence}, the family $\{\nabla p_\gamma\}$ is uniformly bounded in 
$L^2(Q_T)$.  Moreover, it follows that $\nabla p_\gamma \rightharpoonup \nabla p_\infty$ weakly in $L^2(Q_T)$. 
\end{theorem}

\begin{proof}
1. First, we recall the Benilan-Aronson estimates given in Proposition \ref{prop: Benilan-aronson estimates}, that is, 
\begin{align*}
\Theta(u_\gamma) \,  \Delta p_\gamma + \frac{u_{\gamma}}{\gamma^\alpha} \, \Phi(p_\gamma) \geq - \frac{r_{\Phi}}{\gamma^\alpha} \, u_\gamma \, \sigma, 
\end{align*}
where we write shortly $\sigma$ as $\sigma= \sigma(- \frac{r_{\Phi}}{\gamma^{\alpha-1}}t)$, and also that we have assumed, without loss of generality, $\Theta^\prime(\cdot)> 0$.  
Then, we multiply the estimate above by $\Lambda(u_{\gamma})$, where 
$$
   \text{$\Lambda(z):= \frac{z^\gamma}{\Theta(z)}$ when $z \neq 0$, $\Lambda(0)= 0$},
$$
and integrating in $\Omega$ yields
\begin{align*}
\int_{\Omega} p_{\gamma}  \, \Delta p_{\gamma}  \, dx
+ \int_{\Omega} \frac{1}{\gamma^\alpha} u_{\gamma} \, \Lambda(u_{\gamma}) \, \Phi(p_{\gamma}) \, dx \geq - \frac{r_{\Phi}}{\gamma^\alpha} \int_{\Omega} u_{\gamma} \, \Lambda(u_{\gamma}) \, \sigma \, dx. 
\end{align*}
\no In what follows, we integrate by parts the first integral in the right hand side of the inequality above to obtain that
\begin{equation}
\label{ESTGRADPRESSA}
 \int_{\Omega} |\nabla p_{\gamma}|^2 dx \leq  \int_{\Gamma} f_{\gamma} \, \Lambda(u_{\gamma}) \, dr +  \int_{\Omega} \frac{1}{\gamma^\alpha} u_{\gamma} \Lambda(u_{\gamma}) \Phi(p_{\gamma}) \, dx 
 + \frac{r_{\Phi}}{\gamma^\alpha} \int_{\Omega} u_{\gamma} \Lambda(u_{\gamma}) \sigma \, dx. 
\end{equation}

\medskip
2. Now, let us intregrate equation \eqref{ESTGRADPRESSA} from $0$ to $T$, that is, 
\begin{align*}
\iint_{Q_T} |\nabla p_{\gamma}|^2 \, dx dt  \leq C  \Big( \int_{\Gamma_T} |f_{\gamma}| \,d \mathcal{H}^{d} +
 \frac{1}{\gamma^\alpha} \Phi(0)  ||u_{\gamma}||_{L^1(Q_T)} + \frac{r_{\Phi}}{\gamma^\alpha} ||u_{\gamma}||_{L^1(Q_T)} \int_0^T \sigma \, dt \Big),
\end{align*}
where the positive constant $C$ does not depend on $\gamma> 1$. 
Then, we must treat the integral $\int_0^T \sigma dt$ as an improper integral. Indeed, we have 
\begin{align*}
\int_0^T \sigma \Big (- \frac{r_{\Phi} t}{\gamma^{\alpha -1}} \Big) dt &= \displaystyle \lim_{\varepsilon \ra 0} \int_\varepsilon^T  \Big (- \frac{r_{\Phi} t}{\gamma^{\alpha -1}} \Big) dt  \\
&= \gamma^{\alpha-1} \log \Big ( 1- e^{- \frac{r_\Phi T}{\gamma^{\alpha -1}}}  \Big ) -   \gamma^{\alpha-1} \displaystyle \lim_{\varepsilon \ra 0}  \log \Big ( 1- e^{- \frac{r_\Phi \varepsilon}{\gamma^{\alpha -1}}}  \Big ) \\
& = \gamma^{\alpha-1} \log \Big ( 1- e^{- \frac{r_\Phi T}{\gamma^{\alpha -1}}}  \Big ).
\end{align*} 
Collecting all  the above estimates results, it follows that
$$
\begin{aligned}
\iint_{Q_T} |\nabla p_{\gamma}|^2 \, dx dt &\leq C \Big( \int_{\Gamma_T} |f_\gamma| \, d\mathcal{H}^d 
\\[5pt]
&+ \frac{1}{\gamma^\alpha} \Phi(0) ||u_{\gamma}||_{L^1(Q_T)} 
+ \frac{r_{\Phi}}{\gamma} ||u_{\gamma}||_{L^1(Q_T)} \log \big ( 1 - e^{- \frac{r_{\Phi T}}{\gamma^{\alpha-1}}}\big )\Big),
\end{aligned}
$$
from which we obtain
\begin{equation}
\label{GrapPressure}
\sup_{\gamma> 1} \iint_{Q_T} |\nabla p_{\gamma}|^2 \, dxdt < \infty.
\end{equation}

\medskip
3. Finally, due to the weak compactness in $L^2(Q_T)$, there exists a vector field $\mathbf{w} \in L^2(Q_T)$ such that
\begin{align*}
\nabla p_\gamma \rightharpoonup \mathbf{w} \quad \text{weakly in $L^2(\Omega_T)$}  \quad \text{ as } \gamma \ra \infty.
\end{align*}
Since we already have that $p_\gamma$ converges weakly 
as $\gamma \ra \infty$ to $p_\infty \in BV(Q_T)$, it follows  from the uniqueness of the limit that $\mathbf{w}= \nabla p_\infty$. 
\end{proof}

\subsection{Continuity in time and the initial strong trace.}

\begin{proposition} 
\label{prop: continuity in time t=0}
Let the assumptions of Theorem \ref{thm: main theorem} to hold.  Then, 
the weak limit $u_\infty$ given by Theorem \ref{them: theorem convergence} is such that, $u_\infty \in C([0,T]; L^1(\Omega)) \cap BV(Q_T)$ and satisfies the following limit
\begin{align} \label{eq: continuity in t=0 trace eq}
\lim_{t \to 0^+} \int_{\Omega} |u_\infty(t) - u_0| \, dx= 0.
\end{align}
\end{proposition}

\begin{proof}
1. First,  from Theorem \ref{them: theorem convergence}, Theorem \ref{thm: l2 bounds on the gradient of the pressure} we know, 
respectvely that, $u_\infty \in L^1(Q_T)$, and $\nabla p_\infty \in L^2(Q_T)$.   
In order to prove that, $u_\infty \in C([0,T]; L^1(\Omega))$, it is enough to show that $u_{\infty}$ is uniformly 
continuous in time $t$ with values in $L^1(Q_T)$. Therefore, let $I$ be a set of full measure defined by 
$$
   I:= \Big\{ t \in [0,T] / \int_\Omega | \nabla p_\infty(t)|^2 \, dx< \infty  \quad {\rm and} \quad \|f(t)\|_{L^\infty(\Gamma)}< \infty \Big\}. 
$$ 
Take $t_1,t_2 \in I$ such that, $0< t_1 < t_2 \leq T$, and let   
$\varphi \in C^\infty_c(\RR^d)$, $\zeta \in C^\infty_c(\RR)$ be test functions satisfying,
$0 < \zeta(t), \varphi(x) <1$, also $\zeta$ is not a decreasing function on $[0,T]$. 
Therefore, due to Corollary \ref{corollary: pointwise positive limit time derivatives} we have that 
$u_{\infty}$ is not decreasing, and we may write
\begin{equation}
\label{cont_time}
\begin{aligned}
\int_{\Omega} &| u_{\infty}(t_2) \zeta(t_2) - u_{\infty}(t_1) \zeta(t_1)| \, \varphi \, dx 
= \lim_{\gamma \to \infty} \int_{\Omega} \big( u_{\gamma}(t_2) \zeta(t_2) - u_{\gamma}(t_1) \zeta(t_1) \big) \, \varphi \, dx 
\\[5pt] 
& =  \lim_{\gamma \to \infty} \int_{t_1}^{t_2} \!\! \int_{\Omega} \Big( \big( \text{div}\big(\Theta(u_\gamma) \, \nabla p_{\gamma}\big)
+ \frac{1}{\gamma^\alpha} u_\gamma \, \Phi(p_\gamma) \big) \, \varphi  \, \zeta + u_\gamma \zeta^\prime \varphi \Big) \, dx dt 
\\[5pt]
& =   \lim_{\gamma \to \infty} \int_{t_1}^{t_2} \Big( \int_\Gamma f_\gamma \, \varphi \, \zeta \, dr - \int_{\Omega} \zeta \,  \Theta(u_\gamma) \, \nabla p_{\gamma} \cdot \nabla \varphi \, dx 
\\[5pt]
& \qquad \qquad + \int_\Omega \big( \frac{1}{\gamma^\alpha} \, \Phi(p_\gamma)  \, \zeta +  \zeta^\prime \big) \, u_\gamma \, \varphi \, dx \Big) \, dt  
\\[5pt]
& =  \int_{t_1}^{t_2} \Big( \int_\Gamma f \, \varphi  \, \zeta \, dr - \int_{\Omega}  \zeta \, \Theta(u_\infty) \, \nabla p_{\infty} \cdot \nabla \varphi \, dx 
+ \int_{\Omega} u_\infty \, \zeta^\prime \, \varphi \, dx  \Big) \, dt
\\[5pt]
& \leq C \,\big( 1 +  ||\nabla \varphi||_{L^2(\RR^d)}  + {\rm Lip}(\zeta) \big) \, (t_2 - t_1), 
\end{aligned}
\end{equation}
where $C$ is a positive constant. Thus taking sequences of test functions $\{\varphi_n\}$, $\{\zeta_n\}$ that converges to 1, respectively, in $\Omega$
and $[0,T]$, the last inequality yields the desired conclusion that $u_{\infty} \in C([0,T]; L^1(\Omega))$.

\medskip
2. Now, let us identify the trace of $u_\infty$ at $t= 0$. 
Similiar to item (1) above, we may write
$$
\begin{aligned}
\int_{\Omega} &u_{\gamma}(t) \, \zeta(t)  \, \varphi \, dx - \int_\Omega  u_{0\gamma} \, \zeta(0) \, \varphi \, dx 
\\[5pt]
& =  \int_{0}^{t} \Big( \int_\Gamma f_\gamma \, \varphi \, \zeta \, dr - \int_{\Omega} \zeta \,  \Theta(u_\gamma) \, \nabla p_{\gamma} \cdot \nabla \varphi \, dx 
\\[5pt]
& \qquad \qquad + \int_\Omega \big( \frac{1}{\gamma^\alpha} \, \Phi(p_\gamma)  \, \zeta +  \zeta^\prime \big) \, u_\gamma \, \varphi \, dx \Big) \, dt.
\end{aligned}
$$
Then, passing to the limit as $\gamma \to \infty$, we obtain 
$$
\begin{aligned}
\int_{\Omega} &u_{\infty}(t) \, \zeta(t)  \, \varphi \, dx - \int_\Omega  u_{0} \, \zeta(0) \, \varphi \, dx 
\\[5pt]
& =  \int_{0}^{t} \Big( \int_\Gamma f \, \varphi  \, \zeta \, dr - \int_{\Omega}  \zeta \, \Theta(u_\infty) \, \nabla p_{\infty} \cdot \nabla \varphi \, dx 
+ \int_{\Omega} u_\infty \, \zeta^\prime \, \varphi \, dx  \Big) \, dt. 
\end{aligned}
$$
Finally, letting first $t \to 0$ in the above equality, and then $\zeta, \varphi \to 1$, we conclude that 
$$
   u_\infty(0)= u_0 \quad \text{in $L^1(\Omega)$},
$$
which ends the proof.
\end{proof}

\subsection{The incompressibility condition.}

\begin{theorem} \label{them: incompressibility solution}
Let $(u_\infty, p_\infty)$ be given by Theorem \ref{them: theorem convergence}. Then, the pair $(u_\infty, p_\infty)$ satisfies the incompressibility condition  \eqref{V21}.
\end{theorem}

\begin{proof}
1. Firstly let us recall that, $u_\infty \in C([0,T], L^1(\Omega))$,
the Hele-Shaw graph property, that is $p_\infty (u_\infty -1)= 0$, and $0 \leq u_\infty \leq 1$ 
almost everywhere in $Q_T$. Moreover, we have $\nabla p_\infty \in L^2(Q_T)$. In order to show the 
incompressibility condition \eqref{V21}, we may assume
that $p_\infty> 0$ almost everywhere in $Q_T$. Indeed, if $p_\infty \equiv 0$ in a set of positive Lebesgue measure,
then $\nabla p_\infty= 0$ a.e. and we are done on it. 

\medskip
2. Now, since $p_\infty> 0$ due to the Hele-Shaw graph property, 
for sufficiently large $\gamma_0 > 1$ one has for any $\gamma\geq \gamma_0$ that,
for almost all $(t,x) \in Q_T$, 
$u_{\gamma}$ is extremely  close to one,
$g(u_{\gamma}) \geq g_0$ for some $g_0>0$, and also we have 
$p_\gamma \geq p_0$ for some $p_0>0$.  
Then, multiplying the equation
\begin{align*}
\partial_t u_{\gamma} = \text{div} (\Theta(u_{\gamma}) \nabla p_{\gamma}) + \frac{1}{\gamma^\alpha} u_{\gamma} \Phi(p_\gamma)
\end{align*}
by $\gamma \, u_{\gamma}^{\gamma-1}$ yields
$$
\begin{aligned}
\partial_t p_{\gamma} &= \gamma \, u_{\gamma}^{\gamma-1} \, \text{div} \Big ( g(u_{\gamma}) \frac{\nabla p_{\gamma}}{h(u_{\gamma})}\Big ) + \frac{1}{\gamma^{\alpha-1}} p_{\gamma} \Phi(p_\gamma) 
\\[5pt]
& =  g'(u_{\gamma}) \frac{|\nabla p_{\gamma}|^2}{h(u_{\gamma})} 
+ \gamma \, u_{\gamma}^{\gamma-1} \, g(u_{\gamma}) \, \text{div} \Big ( \frac{\nabla p_{\gamma}}{h(u_{\gamma})} \Big ) 
+ \frac{1}{\gamma^{\alpha-1}} \, p_{\gamma}  \Phi(p_\gamma).
\end{aligned}
$$
To follow, let $\zeta \in C_c^\infty(0,T)$, $\xi \in C_c^\infty(\Omega)$ be non-negative test functions.  
Then we multiply the equation above by $\zeta(t)$ $\xi(x)$ and integrate in $Q_T$ to obtain 
\begin{equation}
\label{INCOMP}
\begin{aligned}
\int_0^T \zeta(t) \int_\Omega \frac{\nabla p_\gamma}{h(u_\gamma)} \cdot \nabla \xi(x) \, dxdt
=& -\frac{1}{\gamma} \iint_{Q_T} \frac{u_\gamma}{p_\gamma \, g(u_\gamma)} \, \partial_t p_\gamma \, \zeta(t) \xi(x) \, dxdt
\\[5pt]
& + \frac{1}{\gamma} \iint_{Q_T} \frac{u_\gamma}{p_\gamma \, g(u_\gamma)} \, \frac{g^\prime(u_\gamma)}{h(u_\gamma)} \, |\nabla p_\gamma|^2 \, \zeta(t) \xi(x) \, dxdt
\\[5pt]
& +\frac{1}{\gamma^\alpha} \iint_{Q_T} \frac{u_\gamma}{g(u_\gamma)} \, \Phi(p_\gamma) \, \zeta(t) \xi(x) \, dxdt. 
\end{aligned}
\end{equation}
For each $\gamma \geq \gamma_0$ we obtain the following two estimates from \eqref{INCOMP}: 

\medskip
\underline{First estimate.}
\begin{equation}
\label{INCOMPBelow}
\begin{aligned}
\int_0^T \zeta(t) \int_\Omega \frac{\nabla p_\gamma}{h(u_\gamma)} \cdot \nabla \xi(x) \, dxdt
\geq & -\frac{\sqrt[\gamma]{p_M}}{\gamma} \frac{1}{p_0 \, g_0} \iint_{Q_T}  \partial_t p_\gamma \, \zeta(t) \xi(x) \, dxdt
\\[5pt]
& +\frac{1}{\gamma^\alpha} \iint_{Q_T} \frac{u_\gamma}{g(u_\gamma)} \, \Phi(p_\gamma) \, \zeta(t) \xi(x) \, dxdt,
\end{aligned}
\end{equation}
where we have used \eqref{CondGPrime}. 

\medskip
\underline{Second estimate.}
\begin{equation}
\label{INCOMPAbove}
\begin{aligned}
\int_0^T \zeta(t) \int_\Omega \frac{\nabla p_\gamma}{h(u_\gamma)} \cdot \nabla \xi(x) \, dxdt
\leq & -\frac{1}{\gamma}\;  \frac{1}{2 \, p_M \, g(\sqrt[\gamma]{p_M})} \iint_{Q_T}  \partial_t p_\gamma \, \zeta(t) \xi(x) \, dxdt
\\[5pt]
& + \frac{1}{\gamma} \; \frac{C}{p_0 \, g_0 \, h_0} \iint_{Q_T} |\nabla p_\gamma|^2 \, dxdt
\\[5pt]
& +\frac{1}{\gamma^\alpha} \iint_{Q_T} \frac{u_\gamma}{g(u_\gamma)} \, \Phi(p_\gamma) \, \zeta(t) \xi(x) \, dxdt, 
\end{aligned}
\end{equation}
where the positive constant $C$ does not depend on $\gamma$. 

\medskip
3. Finally, we pass to the limit inferior in the first estimate \eqref{INCOMPBelow} and similarly 
to the limit superior in the second estimate \eqref{INCOMPAbove}, then we obtain 
$$
  0 \leq \liminf_{\gamma \to \infty} \int_0^T \zeta(t) \int_\Omega \frac{\nabla p_\gamma}{h(u_\gamma)} \cdot \nabla \xi(x) \, dxdt
  \leq \limsup_{\gamma \to \infty} \int_0^T \zeta(t) \int_\Omega \frac{\nabla p_\gamma}{h(u_\gamma)} \cdot \nabla \xi(x) \, dxdt \leq 0, 
$$
where we have used \eqref{GrapPressure} and 
$$
   \lim_{\gamma \to \infty} \partial_t p_\gamma= \partial_t p_\infty \quad  \text{weakly in $\mM^1(Q_T)$}. 
$$
Consequently, we have for each non-negative test functions $\zeta \in C_c^\infty(0,T)$, $\xi \in C_c^\infty(\Omega)$, 
$$
   \int_0^T \zeta(t) \int_\Omega \frac{\nabla p_\infty}{h(u_\infty)} \cdot \nabla \xi(x) \, dxdt= 0, 
$$
from which follows the incompressibility condition  \eqref{V21}. Indeed, let $\Pi$ be a distribution such that, 
$\Pi(\chi)= 0$ for any non-negative test function $\chi$, then $\Pi(\phi)= 0$ for any test function $\phi$.
It is enough to take $\phi= \phi^+ - \phi^-$, where $\phi^+, \phi^-$ are respectively the positive and negative
parts of $\phi$.  

\end{proof}

\section{Proof of Main Theorem}

The main issue of this section is to collect all the results of the previous section 
and establish the proof of the Theorem \ref{thm: main theorem} (Main Theorem). 
In particular, we perform a limit transition as $\gamma \to \infty$, (stiff limit), and show existence 
of weak solutions for the Buckley-Leverett System in the sense of Definition \ref{DEFSOL}. 

\begin{proof}[\textbf{Proof of Theorem \ref{thm: main theorem}}]
1. First, let us show item 1. The existence and uniqueness of strong solutions $(u_\gamma, p_\gamma)$
for the compressible Buckley-Leverett System \eqref{stiffpressure}, \eqref{eq: the compressible BL system} and \eqref{BCIC}
follows from Theorem \ref{thm: existence uniq strong solution of the compressible BL}. The strong convergence 
$(u_\gamma, p_\gamma) \to (u_\infty, p_\infty)$ in $L^1(Q_T)$ as 
$\gamma \to \infty$, with $u_\infty, p_\infty \in BV(Q_T)$ is due to Theorem \ref{them: theorem convergence}. 
The regularity $u_\infty \in C([0,T]; L^1(\Omega))$ is given by 
Proposition \ref{prop: continuity in time t=0}. 

\medskip
2.  The maximum principle $0 \leq u_\infty \leq 1$, $0 \leq p_\infty \leq p_M$ a.e. in $Q_T$, and also the Hele-Shaw graph property 
$p_\infty (u_\infty - 1)= 0$ a.e. in $Q_T$, are obtained in Proposition \ref{prop: Hele Shaw limit}. 

\medskip
3. The weak convergence $\nabla p_\gamma \rightharpoonup \nabla p_\infty$ in $L^2(Q_T)$ as $\gamma \to \infty$ is given by 
Theorem \ref{thm: l2 bounds on the gradient of the pressure}. 

\medskip
4. Now, let us show that, the pair $(u_\infty, p_\infty)$ solves the Buckley-Leverett System in the sense of Definition \ref{DEFSOL}.
To this end, let $\phi \in C_{c}^{\infty}((0,T)\times \mathbb{R}^{d})$ be any test function, and from the existence of $(u_\gamma, p_\gamma)$
given by item 1, we multiply 
equation \eqref{eq: the compressible BL system} by $\phi(t,x)$ integrate in $Q_T$ and after integrating by parts
we have
\begin{align}
 \label{DGV21Gamma}
\iint_{Q_T}& \big(u_\gamma \;\phi _{t} + \,\Theta(u_\gamma) \; \nabla p_\gamma \cdot \nabla \phi \big)\ d{x} \,dt = 
\frac{1}{\gamma^\alpha} \, \iint_{Q_T}  u_\gamma \, \Phi(p_\gamma) \, dx dt
+ \int_{\Gamma_T} f_\gamma \, \phi \, d\mathcal{H}^d. 
\end{align}
Also from item 1 we known that $(u_\gamma, p_\gamma) \to (u_\infty, p_\infty)$ strongly in $L^1(Q_T)$, then passing to a subsequence,
if necessary, it follows that $u_\gamma$, $p_\gamma$ converge almost everywhere to $u_\infty$, $p_\infty$ respectively. 
Therefore, the Dominated Convergence Theorem allow us to pass to the limit as $\gamma \to \infty$ in \eqref{DGV21Gamma} obtaining
\eqref{DGV21}, that is, 
$$
    \iint_{Q_T} \big(u_\infty \;\phi _{t} + \,\Theta(u_\infty) \; \nabla p_\infty \cdot \nabla \phi \big)\ d{x} \,dt = \int_{\Gamma_T} f \, \phi \, d\mathcal{H}^d,  
$$
where we have used \eqref{DATABOUNDARY}. Moreover, we have from Proposition \ref{prop: continuity in time t=0} that the inital data is attained 
in the $L^1-$strong sense, that is to say, 
$$
  \lim_{t \to 0^+} \int_{\Omega} |u_\infty(t) - u_0| \, dx= 0.
$$
It remais to show the incompressibility condition  \eqref{V21}, which follows from Theorem \ref{them: incompressibility solution},
that is, for a.a. $t \in (0,T)$ and any text funtion $\xi \in C^\infty_c(\Omega)$, we have 
$$
  \int_{\Omega} \frac{\nabla p_\infty(t)}{h(u_\infty(t))} \cdot  \nabla \xi \, dx= 0. 
$$

\medskip
5. Finally,  the important qualitative properties for the limit, that is, 
$$
   \partial_t u_\infty, \partial_t p_\infty \geq 0 \quad \text{in $\mathcal{D}^\prime(Q_T)$}
$$ 
follows from Corollary \ref{corollary: pointwise positive limit time derivatives}.

\end{proof}

\section*{Data availability statement}
Data sharing is not applicable to this article as no data sets were generated or analysed during the current study.

 \section*{Conflict of Interest}
The author Andr\'{e} de Oliveira Gomes acknowledges and thanks the financial support from the FAPESP
grant number 2018/06531-1 at the University of Campinas (UNICAMP), SP-Brazil.
The author Wladimir Neves has received research grants from CNPq
through the grants  308064/2019-4, 406460/2023-0, and also by FAPERJ 
(Cientista do Nosso Estado) through the grant E-26/201.139/2021.  

\bibliography{Edited_Bibliograph}
\bibliographystyle{plain}

\end{document}